\documentclass[12pt, reqno]{amsart}
\usepackage{tikz}

\usepackage{amsmath,amssymb,enumerate}

\usepackage{epsfig,fancyhdr,color}%,showkeys,amsmidx

\usepackage{amssymb}
\usepackage{amsmath,amsthm}
\usepackage{bbm}
\usepackage{latexsym}
\usepackage{amscd}
\usepackage{psfrag}
\usepackage{graphicx}
\usepackage[latin1]{inputenc}
\usepackage[all]{xy}

\usepackage[flushleft]{threeparttable}
\usepackage{subcaption}
\usepackage{epsfig,fancyhdr,color}%,howkeys,amsmidx
\usepackage{amssymb}
\usepackage{amsmath,amsthm}

\usepackage{mathtools}

\usepackage[pdfpagemode={UseOutlines},bookmarks=true,bookmarksopen=true, bookmarksopenlevel=0,bookmarksnumbered=true,hypertexnames=false, colorlinks,linkcolor={blue},citecolor={blue},urlcolor={blue}, pdfstartview={FitV},unicode,breaklinks=true,backref=page]{hyperref}

\definecolor{NoteColor}{rgb}{1,0,0}

\newtheorem{theorem}{\rm\bf Theorem}[section]
\newtheorem{proposition}[theorem]{\rm\bf {Proposition}}
\newtheorem{lemma}[theorem]{\rm\bf Lemma}
\newtheorem{corollary}[theorem]{\rm\bf Corollary}
\newtheorem*{theorem 1}{\rm\bf {proposition} 1}
\newtheorem*{theorem 2}{\rm\bf {proposition} 2}

\theoremstyle{definition}

\theoremstyle{remark}
\newtheorem{remark}[theorem]{\rm\bf Remark}

\def\interieur#1{\mathord{\mathop{\kern 0pt #1}\limits^\circ}}

\def\hyperp{{\rm I}\kern-.3ex{\rm H}}

\def\hyperl{\mathbb H}

\usepackage{mathabx,epsfig}
\def\acts{\mathrel{\reflectbox{$\righttoleftarrow$}}}

\begin{document}

\title[Super Hyperbolic Law of Cosines]{Super Hyperbolic Law of Cosines: same formula with
different content}

\author{Robert Penner}
\address {\hskip -2.5ex Institut des Hautes \'{E}tudes Scientifiques\\
35 route des Chartres\\
Le Bois Marie\\
91440 Bures-sur-Yvette\\
France\\
{\rm and}~Mathematics Department,
UCLA\\
Los Angeles, Ca 90095\\USa}
\email{rpenner{\char'100}ihes.fr}
\


 \date{\today}
 
 \thanks{Thanks to Thibault Damour and Igor Frenkel for helpful comments and likewise to Anton Zeitlin as well as for earlier induction into the world of super mathematics.}
 

% ---------------------------------------

\begin{abstract}
We derive the Laws of Cosines and Sines in the super hyperbolic plane using Minkowski supergeometry and find the identical formulae to the classical case, but remarkably involving different expressions for the trigonometric functions, which include substantial fermionic corrections.  In further analogy to the classical case, we apply these results to show that two parallel supergeodesics which are not ultraparallel admit a unique common orthogonal supergeodesic, and we briefly describe aspects of elementary supernumber theory, leading to a prospective analogue of the Gauss product of quadratic forms. 
\end{abstract}

\maketitle

\section*{Introduction}

\noindent It is a pleasure and and an honor to participate in this volume celebrating Norbert A'Campo on the occasion of his 80th birthday.  We have been friends for roughly one quarter of his life, approximately one third of mine, enthusiastically introduced to one another by our common lifelong friend Athanase Papadopoulos.  Athanase was certain we would be sympatico, and he was correct.

\medskip

\noindent There are many ways to do mathematics, from abstracting ponderously huge machinery to computing explicit examples
and everything in between. I have learned that the A'Campo way is something quite special: To appreciate a field of flowers, one can ascertain the deeper beauty and structure of a single bloom, its profound reality, and cultivate a preternatural intimacy with it.  Nothing fancy, nothing grand, but rather a humble and natural communion that nevertheless brings with it the deepest comprehension
of the entire meadow as reflected in one lone blossom. 

\medskip

\noindent
In a more mathematical rendition, Athanase paraphrases the great Lobachevsky as follows:

\medskip

{\sl In order to understand a geometry, both locally and globally, it is enough to understand its triangles.}

\medskip

\noindent  This informs the basic purview of my birthday present here to Norbert, to answer the question: 

\medskip

{\sl What are the Laws of Cosines and Sines for triangles in hyperbolic superspace?}

\medskip

\noindent In order to explain,
let ${\mathbb K}_*={\mathbb K}_*[0]\otimes{\mathbb K}_*[1]$ be the ${\mathbb Z}$/2-graded algebra over the field ${\mathbb K}={\mathbb R}$ or ${\mathbb C}$ with one central generator
$1\in {\mathbb K}\subseteq {\mathbb K}_*[0]$ of degree zero and countably infinitely many anti-commuting generators in ${\mathbb K}_*[1]$ of degree one.  $a_\#$ is called the {\it body} of the {\it supernumber} $a\in {\mathbb K}_*$,
and $a$ is said to be {\it even} or {\it odd}, respectively, if it lies in ${\mathbb K}_*[0]$  or ${\mathbb K}_*[1]$.

\medskip

\noindent One can define Riemannian supermanifolds and super Riemann surfaces in the natural way modeled on affine spaces constructed from $K_*$, as in \cite {dewitt} and see also \cite{Rogers} in general,  and as in \cite{CraneRabin,PZ,Witten} for Riemann surfaces in particular.
This is discussed in the next section in detail sufficient for our purposes here.  While no particular
mathematical justification for such formalism is necessary, for one can study a flower simply for its own reward especially according to the A'Campo way, it is worth pointing out that supermanifolds are part and parcel of the Standard Model of high energy physics, roughly with odd variables corresponding to the fermions comprising matter, and even variables to the bosons formalizing interactions between them.

\medskip

\noindent A fundamental example  from physics is the super upper half plane  \cite{CraneRabin}
$${\mathcal U}_*=\{ (z,\theta)\in{\mathbb C}^{1|1}:{\mathcal Im}~z={1\over 2}(z-\bar{z})>0\},$$
where the overline denotes complex conjugation, with its supermetric 
$ds=({\mathcal Im}~z+{1\over 2}\theta\bar\theta)^{-1}\vert dz+\theta d\theta\vert$
invariant under the action 
$$\begin{aligned}
z&\mapsto {{az+b}\over{cz+d}}+\theta{{\gamma z+\delta}\over{(cz+d)^2}}\\
\theta&\mapsto{{\gamma z+d}\over{cz+d}}+{{\theta(1+{1\over2}\delta\gamma})\over{cz+d}}\\
\end{aligned}$$
of the so-called M\"obius supergroup $\rm SPL(2,{\mathbb R})$, where $a,b,c,d$ are even satisfying $ad-bc=1$, and $\gamma,\delta$ are
odd.  In particular, elements of the M\"obius supergroup have three even and two odd degrees of freedom.

\medskip 

\noindent The body of ${\rm SPL}(2,{\mathbb R})\acts{\mathcal U}_*$ is evidently the classical action
of the M\"obius group on the standard upper half plane model for the hyperbolic plane, with its superabsolute the ``one-superpoint compactification" ${\mathbb S}^{1|1}$ of the real superaxis ${\mathbb R}^{1|1}$, but this is not the model for the action of orientation-preserving isometries of the super hyperbolic plane that we shall adopt here.  
Rather, we shall work with the orthosymplectic group ${\rm OSp}(1|2)\approx {\rm SPL}(2,{\mathbb R})$ \cite{Kac}, recalled in $\S$\ref{osp}, acting on
super Minkowski space ${\mathbb R}^{2,1|2}$, and in particular on its upper sheet $\hyperp$ of the unit hyperboloid of two sheets, as explained in $\S$\ref{smsc}, which gives our equivariantly isometric model ${\rm OSp}(1|2)\acts\hyperp$ for ${\rm SPL}(2,{\mathbb R})\acts{\mathcal U}_*$.

\medskip

\noindent The body of ${\rm OSp}(1|2)\acts{\mathbb R}^{2,1|1}$ is the usual action ${\rm SO}_+(2,1)\acts{\mathbb R}^{2,1}$ on Minkowski 3-space, with its quadratic form
$z^2-x^2-y^2$ in the standard coordinates, under the component ${\rm SO}_+(2,1)$ of the identity
in ${\rm SO}(2,1)$.  Classically, the upper sheet of the hyperboloid of two sheets inherits a Riemannian metric from the ambient quadratic form, giving a model for the hyperbolic plane.  Moreover, the unit hyperboloid of one sheet parametrizes the oriented hyperbolic geodesics, and the open positive light-cone of isotropic vectors at positive height parametrizes the space of all horocycles; cf.~\cite{pennerbook}.  There are analogously defined respective conics
$\hyperl$ and ${\mathcal L}^+$ in Minkowski superspace.

\medskip

\noindent There are numerous reasons for working in the super Minkowski model.  First of all, hyperbolic supergeometry is the same but different from the usual hyperbolic geometry, and some of these differences are more readily manifest in the Minkowski model.  As is well-known \cite{Manin} and clear in either model, orientation-preserving isometries do not act transitively
on triples of ideal points, and orbits of such triples have a single odd modulus, which is the genesis of the odd coordinates on super Teichm\"uller space \cite{PZ}, whose analysis is based on the Minkowski model and follows the general approach of decorated Teichm\"uller theory \cite{pennerbook}.

\medskip

\noindent A less well-known difference, again from \cite{PZ}, is that ${\rm OSp}(1|2)$ does not act transitively on the super positive light-cone; there is again one odd modulus, and the ${\rm OSp}(1|2)$-orbit corresponding to the vanishing of this modulus is the so-called special light-cone ${\mathcal L}_0\subseteq{\mathcal L}^+$ of \cite{PZ}.  It is the action
${\rm OSp}(1|2)\acts{\mathcal L}_0$ whose body captures the usual action of the M\"obius supergroup on the hyperbolic absolute in the sense that there is an equivariant map from the former to the superabsolute, as for the bodies in the classical case.

\medskip

\noindent  One of the most stark differences is that not every supergeodesic ray in ${\mathcal U}_*$ is asymptotic to a point in the superabsolute!
In fact, supergeodesic rays in \hyperp ~are always asymptotic to a ray in ${\mathcal L}^+$, but in ${\mathcal U}_*$ are only asymptotic to a point in the superabsolute if the corresponding supergeodesic ray in $\hyperp$ is asymptotic to a ray in the super positive light-cone ${\mathcal L}^+$ that happens to lie in the special light-cone
${\mathcal L}_0$.  It is still true, however, that two distinct points in the superabsolute determine a unique so-called special supergeodesic which is asymptotic to them, but supergeodesics are generically not special.  The assertions of this paragraph
can be derived from results in \cite{PZ} together with those of this paper, but these matters, together with the further fact that
special supergeodesics are precisely those that are dual to points in the super hyperboloid $\hyperl$ of one sheet, will be taken up elsewhere.

\medskip

\noindent Having thus sketched some of the ways in which classical intuition can be misleading in hyperbolic supergeometry, we turn finally to the scope of this paper.  The first basic question that has motivated this work is what relations hold between the angles
and pairwise distances for a triple of non-collinar points in $\hyperp$.  Of course, the preliminary issues of what are the natures of collinearity, studied in $\S$4, angles and distances in hyperbolic superspace must be addressed first. 
The answer in $\S$5, in a sense disappointing, is that the {\sl exact same formulae} for the classical hyperbolic Laws of Cosines and Sines hold in hyperbolic superspace.  

\medskip

\noindent However, the possible disappointment fades upon realizing that the cosine and sine of an angle in hyperbolic superspace are wildly different from their classical counterparts, as proved in $\S$6, with substantial odd contributions, called {\sl fermionic corrections} as motivated by physics.  Thus, {\sl the same formulae relate wildly different quantities}, and this is quite remarkable.  

\medskip

\noindent These formulae are our main takeaways, but we also employ these results in $\S$6 to include another super analogue of a classical result:
two parallel supergeodesics which are not ultraparallel, in the sense that they do not share asymptotes in ${\mathcal L}^+$, admit a unique common orthogonal supergeodesic.  Moreover along the way in $\S$3, we briefly consider elementary supernumber theory, and in particular an abelian group structure analogous
to the Gauss groups of quadratic forms.  

\medskip

\noindent The reader may correctly gather that the results of this paper are just a beginning, one flower in the field of supergeometry, hyperbolic and otherwise.  It is finally worth mentioning that much of this paper applies to supergeodesics and supertriangles in hyperbolic superspace of arbitrary dimension.

\section{Super Minkowski space and its quadrics}\label{smsc}

\noindent Let ${\mathbb K}_*={\mathbb K}_*[0]\otimes{\mathbb K}_*[1]$ be the ${\mathbb Z}$/2-graded module over the field ${\mathbb K}={\mathbb R}$ or ${\mathbb C}$ with one central generator
$1\in {\mathbb K}\subseteq {\mathbb K}_*[0]$ of degree zero and countably infinitely many anti-commuting generators $e_1,e_2,\ldots\in {\mathbb K}_*[1]$ of degree one, so that an arbitrary $a\in{\mathbb K}_*$ can be written uniquely as
$$a=a_{\#}+\sum_i a_ie_i+\sum_{i<j} a_{ij} e_i e_j++\sum_{i<j<k} a_{ijk} e_i e_j e_k\cdots,$$
where $a_\#, a_i,a_{ij},a_{ijk}\ldots\in {\mathbb K}$.  $a_\#$ is called the {\it body} of the {\it supernumber} $a\in {\mathbb K}_*$.
One allows only finitely many anti-commuting factors in any product and only finitely many summands in any supernumber, a constraint we shall call {\it regularity}.  (The reason for taking infinitely many anti-commuting variables is that in the full theory one wants the graded version $d(fg)=(df)g-f(dg)$ of the usual Leibnitz rule to hold, and this would be confounded by taking $f$ to be the finite product of all the anti-communting variables and $g$ to be any of them.)

\medskip

\noindent Regularity implies that $a-a_\#$ is always nilpotent, i.e., for each $a\in{\mathbb K}_*$, there is some $n\in{\mathbb Z}_{\geq 0}$ with $(a-a_\#)^{n}=0$.  It follows that if $a_\#\neq 0$, then we may write
$$\begin{aligned}
{1\over a}&= {1\over{a_\#+(a-a_\#)}}
={1\over a_\#}~~ {1\over{1+{{a-a_\#}\over a_\#}}}\\
&={1\over a_\#}\biggl [1-{{a-a_\#}\over a_\#}+\bigg({{a-a_\#}\over a_\#}\biggr)^2-\cdots \pm\biggl({{a-a_\#}\over a_\#}\biggr)^{n-1}
\biggr ],\\
\end{aligned}$$
and hence $a\in{\mathbb K}_*$ is invertible if and only if $a_\#\neq 0$.  It similarly follows from regularity that the zero divisors in 
${\mathbb K}_*$ are given by the ideal generated by ${\mathbb K}_*[1]$.  

\medskip

\noindent One analogously extends real-analytic functions
on ${\mathbb K}$ to ${\mathbb K}_*$ with Taylor series under appropriate restrictions on the body.  For instance for later application, we have the elementary

\begin{lemma} \label{quadratics}
Suppose that $\lambda^3=0$.  If $a$ has nonzero body, then
$${1\over{a+\lambda}}={1\over a}\biggl ( 1-{\lambda\over a}+{{\lambda^2}\over {a^2}}\biggr ),$$
and if $a$ has positive body, then
$$\sqrt{a+\lambda}=\sqrt{a}~\biggl ( 1+{1\over 2}{\lambda\over a} -{1\over 8}{\lambda^2\over a^2}\biggr).$$
\end{lemma}

\medskip

\noindent If $a\in{\mathbb K}_*[0]$, then it is said to be an {\it even} supernumber or  {\it boson}, while if $a\in{\mathbb K}_*[1]$, then it is said to  be an {\it odd} supernumber or {\it fermion}.  An order relation
$\leq$ on ${\mathbb K}$ induces one on $a,b\in{\mathbb K}_*$, where $a\leq b$ if and only if $a_\#\leq b_\#$.

\medskip

\noindent  {\it Affine ${\mathbb K}$ superspace} is defined to be
$${\mathbb K}^{n\vert m} = \{(x_1,x_2,\ldots,x_n~\vert~\theta_1,\theta_2,\ldots ,\theta_m)\in {\mathbb K}_*^{n+m}: x_i\in K_*[0], \theta_j\in K_*[1]\}.$$
One can define n$\vert$m {\it supermanifolds} with charts based on affine superspace 
${\mathbb K}^{n|m}$ in the usual way, and a {\it (Riemannian) supermetric} on an $n\vert m$ supermanifold is defined to be a  positive definite boson-valued quadratic form on each tangent space ${\mathbb K}^{n|m}$ as usual.

\medskip

\noindent A tuple of points ${\bf x}_i\in {\mathbb K}^{n\vert m}$ is said to be {\it linearly independent} (in the super sense) if there is no relation $0=\sum_i a_i{\bf x_i}$
where all the $a_i\in{\mathbb K}_*$ have non-zero body.  In particular for $n=2$, ${\bf x}_1$ and
${\bf x}_2$ are said to be {\it parallel} if there is some $a\in{\mathbb K}_*$ with non-zero body
so that ${\bf x}_1=a{\bf x}_2$.

\medskip

\noindent
The principal example for us here is 
{\it (real) super Minkowski {\rm 2,1$\vert$2} space}
$${\mathbb R}^{2,1|2}=\{ (x_1,x_2,y~\vert~\phi,\psi)\in{\mathbb R}_*^5:x_1,x_2,y\in{\mathbb R}[0]~{\rm and}~ \phi,\psi\in{\mathbb R}[1]\},$$
which supports the boson-valued symmetric bilinear pairing
$$\langle (x_1,x_2,y~\vert~\phi,\psi),(x_1',x_2',y'~\vert~\phi',\psi')\rangle={1\over 2}(x_1x_2'+x_1'x_2)-yy'+\phi\psi'+\phi'\psi$$
with associated quadratic form $x_1x_2-y^2+2\phi\psi$.

\medskip

\noindent The body of ${\mathbb R}^{2,1|2}$ with this inner product is evidently
the classical Minkowski space ${\mathbb R}^{2,1}$ with its (negative) definite restriction to
$$\begin{aligned}
\hyperp'&=\{ {\bf x}=(x_1,x_2,y)\in{\mathbb R}^{2,1}:\langle{\bf x},{\bf x}\rangle=1~{\rm and}~x_1+x_2>0\}\\
\end{aligned}$$
providing a model of the hyperbolic plane, and
we analogously define the {\it super hyperbolic plane} to be
$$\begin{aligned} 
{\hyperp}&=\{ {\bf x}=(x_1,x_2,y~\vert~\phi,\psi)\in{\mathbb R}^{2,1|2}:\langle{\bf x},{\bf x}\rangle=1~{\rm and}~x_1+x_2>0\}\supseteq\hyperp'\\
\end{aligned}$$ 
with its metric likewise induced from the inner product.

\medskip

Continuing by analogy,
let $$\hyperl'=\{ {\bf h} \in{\mathbb R}^{2,1}:\langle{\bf h},{\bf h}\rangle=-1\}$$ denote the hyperboloid of one sheet
and 
$${\hyperl}=\{ {\bf h} \in{\mathbb R}^{2,1|2}:\langle{\bf h},{\bf h}\rangle=-1\}\supseteq\hyperl'$$
its super analogue.  
Let $L$ denote the collection of isotropic vectors in ${\mathbb R}^{2,1}$ with
$$L^+=\{ {\bf u}=(u_1,u_2,v)\in{L}:~u_1+u_2>0\}$$ 
 the {\it (open) positive light cone} (whose points are affine duals to horocycles in $\hyperp'$, as in \cite{pennerbook}), and let 
${\mathcal L}$ denote the collection of isotropic vectors in ${\mathbb R}^{2,1|2}$ with
$${\mathcal L}^+=\{ {\bf u}=(u_1,u_2,v~\vert~\xi,\eta)\in{\mathcal L}:u_1+u_2>0\}.$$
Though we shall not require it here, the {\it special light-cone} discussed in the Introduction is defined by
$$\begin{aligned}
{\mathcal L}_0&=\{{\bf u}=(u_1,u_2,v~\vert~\phi,\psi)\in{\mathcal L}^+:u_1\psi=v\phi\};
\end{aligned}$$
note the equivalence of the conditions $u_1\psi=v\phi$ and $u_2\phi=v\psi$ on ${\mathcal L}$.

\bigskip

\section{The orthosymplectic group OSp(1$\vert$2)}\label{osp}

Here we provide basic information concerning the orthosymplectic group, which is the simplest Lie supergroup, whose body is the classical special linear group ${\rm SL}(2,{\mathbb R})$.  We refer the interested reader to \cite{Kac,Manin}
for further details about general Lie superalgebras and supergroups
and to \cite{PZ} for details about OSp(1$\vert$2)$\approx$SPL(2,${\mathbb R}$). 

\medskip

Given a Lie algebra $\mathfrak{g}$, consider 
the Lie superalgebra $\mathfrak{g}(S)=S\otimes\mathfrak{g}$ for some Grassmann algebra $S$ with its decomposition $S=S[0]\oplus S[1]$
into even and odd elements.
It follows that $\mathfrak{g}(S)$ is both a right and left $S$-module, i.e.,
$s\otimes T=(-1)^{|s||T|}(1\otimes T)(s\otimes 1)$ if $s\in S$ and $T\in\mathfrak{g}$ are homogeneous elements of respective degrees $|s|$ and $|T|$.
This provides a representation of the corresponding Lie superalgebra $\mathfrak{g}(S)$ in the space $S\otimes\mathbb{R}^{m|n}$ from a given representation of $\mathfrak{g}$ in $\mathbb{R}^{m|n}$, and then a representation of the corresponding Lie supergroup $G(S)$ by exponentiating pure even elements from $\mathfrak{g}(S)$ in $S\otimes\mathbb{R}^{m|n}$.  

\medskip

In particular to be entirely explicit about the signs when writing a super matrix representing the action of 
 ${\rm SL}(2,{\mathbb R})(S)$ or $\mathfrak{sl_2}(S)$ as elements of $S\otimes End(\mathbb{R}^{m|n})$ on $S\times\mathbb{R}^{m|n}$, the product
 of two supermatrices from OSp(1$\vert$2) is given by\footnote{The usual (super)matrix multiplication (without the minus signs above) is recovered upon replacing the odd entries in the third row by their negatives. This difference in sign arises from the fact that one typically considers the action of group elements on ${S[0]}^{m}\times {S[1]}^{n}$, which can be identified with the space of even elements in $S\otimes \mathbb{R}^{m|n}$, and the extra minus sign comes from that isomorphism.}

 $$\biggl ( \begin{smallmatrix}
a_1&b_1&\alpha_1\\c_1&d_1&\beta_1\\\gamma_1&\delta_1&f_1\\
\end{smallmatrix}\biggr )
\biggl ( \begin{smallmatrix}
a_2&b_2&\alpha_2\\c_2&d_2&\beta_2\\\gamma_2&\delta_2&f_2\\
\end{smallmatrix}\biggr )=
\biggl (\begin{smallmatrix}
a_1a_2+b_1c_2-\alpha_1\gamma_2&a_1b_2+b_1d_2-\alpha_1\delta_2&a_1\alpha_2+b_1\beta_2+\alpha_1 f_2\\
c_1a_2+d_1c_2-\beta_1\gamma_2&c_1b_2+d_1d_2-\beta_1\delta_2&c_1\alpha_2+d_1\beta_2+\beta_1 f_2\\
\gamma_1 a_2+\delta_1 c_2+\delta_1\gamma_2&\gamma_1 b_2+\delta_1 d_2+f_1 \delta_2&-\gamma_1\alpha_2-\delta_1\beta_2+f_1f_2\\
\end{smallmatrix}\biggr ).$$
The {\it superdeterminant} or {\it Berezinian} of $g=\begin{psmallmatrix}
a&b&\alpha\\c&d&\beta\\\gamma&\delta&f\\
\end{psmallmatrix}\in {\rm OSp}(1|2)$ is given by
\[
{\rm sdet}~g~=~f^{-1}
~\det\left[
\begin{pmatrix}
a&b\\ c&d
\end{pmatrix}
+f^{-1}
\begin{pmatrix}
\alpha\gamma&\alpha\delta\\ \beta\gamma&\beta\delta
\end{pmatrix}
\right]
\]
provided $f$ is invertible.  The analogue of the classical determinant, the Berezinian  \cite{Kac,Manin,Rogers} is characterized by being a  multiplicative homomorphism
satisfying
sdet exp $g$ = exp$(a+d-f)$, i.e., the exponential of the supertrace of $g$, but unlike the classical determinant,
it is only defined for invertible supermatrices.

\medskip

The supergroup ${\rm OSp}(1|2)$ can be faithfully realized as $(2|1)\times (2|1)$ supermatrices $g$ with sdet equal to unity obeying the relation
\begin{eqnarray}\nonumber
g^{st}Jg=J,
\end{eqnarray}
where 
\begin{eqnarray}\nonumber
J=\left( \begin{array}{ccc}
~0 & 1 & ~0 \\
-1 & 0 & ~0 \\
~0 & 0 & -1 \end{array} \right)
\end{eqnarray}
and where the {\it supertranspose} $g^{st}$ of $g$ is given by
\begin{eqnarray}\nonumber
g=\left( \begin{array}{ccc}
a & b & \alpha \\
c & d & \beta \\
\gamma & \delta & f \end{array} \right)\quad
{\rm implies} \quad
g^{st}=\left( \begin{array}{ccc}
~a & ~c & \gamma \\
~b & ~d & \delta \\
-\alpha & -\beta & f \end{array} \right).
\end{eqnarray}
This provides a simple formula
$$
g^{-1}=J^{-1}g^{st}J=\left( \begin{array}{ccc}
~~d & -b & ~~~~\delta \\
-c & ~~a & -\gamma \\
-\beta & ~~\alpha & ~~f \end{array} \right)
$$
for inversion in ${\rm OSp}(1\vert 2)$ as well as 
leading to the system 
\begin{eqnarray}\nonumber
&&\alpha=b\gamma-a\delta, \quad \beta =d\gamma-c\delta, \quad \hskip 2.8ex f=1+\alpha\beta,\nonumber\\
&&\gamma=a\beta-c\alpha, \quad \delta=b\beta-d\alpha, \quad f^{-1}=ad-bc\nonumber
\end{eqnarray}
of constraints on the entries of $g$, which, together with the demand that ${\rm sdet}~g=1$, completely characterize elements of OSp(1$\vert$2).

\medskip

\noindent There is a canonical inclusion ${\rm SL}(2,{\mathbb R})<{\rm OSp}(1|2)$, which extends to 
${\rm SL}(2,{\mathbb R_*[0]})<{\rm OSp}(1|2)$, given by
$(\begin{smallmatrix}a&b\\c&d\\\end{smallmatrix})\mapsto
( \begin{smallmatrix}a&b&0\\c&d&0\\0&0&1\\\end{smallmatrix})$. 
 It is worth emphasizing that it is not the M\"obius group, but rather the full special linear group that appears here, so that a suitable representation of a Fuchsian
group into ${\rm OSp}(1\vert 2)$ provides, upon taking the body, a representation in 
${\rm SL}(2,{\mathbb R})$, or equivalently a spin structure on the underlying Riemann surface, cf. \cite{natanzon,PZ}.  

\medskip

\noindent It is not difficult to prove the following useful 

\begin{lemma}\label{factor}
Any element of {\rm OSp(1$\vert$2)} can be written uniquely in the form
$$\begin{pmatrix}a&b&0\\c&d&0\\0&0&1\\\end{pmatrix}
\begin{psmallmatrix}1-{{\alpha\beta}\over 2}&0&\alpha\\0&1-{{\alpha\beta}\over 2}&\beta\\\beta&-\alpha&1+\alpha\beta\\\end{psmallmatrix}
=
\begin{psmallmatrix}1-{{\alpha\beta}\over 2}&0&a\alpha+b\beta\\0&1-{{\alpha\beta}\over 2}&c\alpha+d\beta\\c\alpha+d\beta&-(a\alpha+b\beta)&1+\alpha\beta\\\end{psmallmatrix}
\begin{pmatrix}a&b&0\\c&d&0\\0&0&1\\\end{pmatrix},$$
for appropriate fermions $\alpha,\beta$ and bosons $a,b,c,d$ with $ad-bc=1$.
\end{lemma}

\noindent As a point of notation for later utility, for any two fermions $\alpha,\beta$, let $$u(\alpha,\beta)=
\begin{psmallmatrix}1-{{\alpha\beta}\over 2}&0&\alpha\\0&1-{{\alpha\beta}\over 2}&\beta\\\beta&-\alpha&1+\alpha\beta\\\end{psmallmatrix}\in {\rm OSp}(1\vert 2).$$

\bigskip

\section{Hyperbolic supergeometry and supernumber theory}
Minowski three space ${\mathbb R}^{2,1}\approx {\mathbb R}^3$, as the space of binary quadratic forms or binary symmetric bilinear forms, is naturally coordinatized by 
$$A=\begin{pmatrix}z+x&y\\y&z-x\end{pmatrix}=\begin{pmatrix}x_1&y\\y&x_2\end{pmatrix}\in{\mathbb R}^{2|1},$$
where the quadratic form is $(u,v)\mapsto \begin{psmallmatrix}u&v\end{psmallmatrix} A \begin{psmallmatrix}u\\v\end{psmallmatrix}=x_1u^2+2yuv+x_2v^2$,
and $g\in {\rm SO}_+(1,2)$ acts via isometry on $A\in{\mathbb R}^{2,1}$ as change of basis via the adjoint
$$g:A\mapsto g^tAg.$$  This action of ${\rm SO}_+(1,2)\approx {\rm PSL}(2,{\mathbb R})$ as the group of isometries of the hyperbolic plane $\hyperp'$ is a fundamental link between hyperbolic geometry and elementary number theory.

\medskip

\noindent Likewise, ${\mathbb R}^{2,1|2}\approx {\mathbb R}^{3|2}$, as a space of quadratic superforms, is naturally coordinatized
by $$A=\begin{pmatrix}~x_1&~~y&\phi\\~~y&~x_2&\psi\\-\phi&-\psi&0\\\end{pmatrix}\in{\mathbb R}^{2,1|2},$$
where $(u,v,\theta)\mapsto\begin{psmallmatrix}u,v,-\theta\end{psmallmatrix} A \begin{psmallmatrix}u\\v\\\theta\end{psmallmatrix}=x_1u^2+2yuv+x_2v^2$ takes the same values as before but with $u,v$ arbitrary bosons, and $g\in{\rm OSp}(1\vert 2)$ acts on $A$  as change of basis again via the adjoint
$$g:A\mapsto g^{st}Ag.$$  One checks using Lemma \ref{factor} that this action is again isometric and hence restricts to an
action on $\hyperp$ itself.  
In fact according to Theorem 1.2 of \cite{PZ}, the mapping
$$
z={{i-y-i\phi\psi}\over x_2}~~{\rm and}~~
\theta={\psi\over x_2}(1+iy)-i\phi
$$
establishes an equivariant isometry $\hyperp\to{\mathcal U}_*$ which conjugates the action
${\rm OSp}(1\vert 2)\acts\hyperp$ to the action SPL(2,${\mathbb R}$)$\acts{\mathcal U}_*$, just as in the classical case \cite{pennerbook}.

\medskip

\noindent 
The body of this action of ${\rm OSp}(1\vert 2)$ on $\hyperp$ is the classical action of orientation-preserving isometries on the hyperbolic plane, and this extension is our analogous action of this Lie supergroup on the super hyperbolic plane.  This is the main point of this section, the rest of which can be skipped for the sequel.
 
\medskip

\noindent  Before turning attention to hyperbolic supergeometry, let us take a moment to
ponder elementary supernumber theory.  Consider the ring ${\mathbb Z}_*$, defined as in the first sentence of Section \ref{smsc} for the ring ${\mathbb K}={\mathbb Z}$ rather than for a field.  Take the usual definition of divisibility, which is
complicated by the plethora of divisors of unity and zero in ${\mathbb Z}_*$,
 in order to define the least common multiple  and greatest common divisor.  These must constitute subsets rather than elements of ${\mathbb Z}_*$, since the ordering on ${\mathbb Z}_*$ induced from ${\mathbb Z}$ is not even a partial ordering, let alone a linear ordering.  Two supernumbers are {\it relatively prime} if their only common divisors are divisors of unity.  

\medskip

\noindent  
As above, $g\in{\rm SL}(2,{\mathbb Z})$ acts by change of basis on
$\begin{psmallmatrix}a&b/2\\b/2&c\end{psmallmatrix}$, whose corresponding
quadratic form  $ax^2+bxy+cy^2$ is called {\it primitive}
provided $a,b,c$ are pairwise relatively prime.
This action evidently
leaves invariant the discriminant $D=b^2-4ac$ and turns out also to preserve
primitivity.
Gauss defined an abelian product on the finite collection
${\mathcal G}(D)$ of ${\rm SL}(2,{\mathbb Z})$-orbits of 
primitive forms of discriminant $D$, which is essentially the ideal 
class group of ${\mathbb Q}(\sqrt{D})$.  
In fact, the Gauss product extends
to an abelian semigroup structure on $\cup {\mathcal G}(D)$, where the union is
over all discriminants with common square-free kernel, a condition interpreted
geometrically in \cite{PGauss}.  

\medskip

\noindent It is natural to wonder whether these considerations might 
extend to supernumbers and supergeometry.  
The sub supergroup ${\rm OSp}(1|2,{\mathbb Z}_*)<{\rm OSp}(1|2)$
whose entries lie in ${\mathbb Z}_*$ acts as above by change of basis on
$\begin{psmallmatrix}a&b/2&\phi\\b/2&c&\psi\\\phi&-\psi&0\\\end{psmallmatrix}$, where $a,b,c\in{\mathbb Z}_*[0]$
are pairwise relatively prime and $\phi,\psi\in{\mathbb Z}_*[1]$, and it leaves invariant the discriminant
$D=b^2-4ac+8\phi\psi$, which differs from the bosonic case even though the quadratic forms agree.  It is natural to conjecture that there is an abelian semigroup structure
on the analogously defined $\cup{\mathcal G}(D)$ with geometric underpinnings similar to \cite{PGauss}.

\medskip

\noindent This putative super version of the Gauss product might be uninteresting, or it might shed light on ideal class numbers.  In any case, the idea of studying supernumber theory seems to warrant further thought.

\bigskip

\section{Supergeodesics}

To begin, we recall from \cite{HPZ} that hyperbolic supergeodesics have the same general parametric form as relativistic geodesics.

\begin{theorem}[Theorem 1.2 of \cite{HPZ}]\label{uvthm}
Geodesics parametrized by arc length on the super hyperboloid of two sheets are described by the equation
$$
{\bf x}={\bf u}\cosh s+{\bf v}\sinh s,
$$
where ${\bf u}\in\hyperp$, ${\bf v}\in\hyperl$ and $\langle{\bf u}, {\bf v}\rangle =0$. 
The asymptotes of the corresponding supergeodesic are given by the rays containing the vectors
$$
\bf e=\bf u+ \bf v, \quad \bf f=\bf u-\bf v\in{\mathcal L}^+
$$
Conversely,  points ${\bf e},{\bf f}\in{\mathcal L}^+$ with $\langle {\bf e},{\bf f}\rangle =2$ uniquely define a geodesic, where ${\bf e}, {\bf f}$ give rise to ${\bf u}={1\over 2}({\bf e}+{\bf f})\in \hyperp$ and ${\bf v}={1\over 2}({\bf e}-{\bf f})\in\hyperl$.   
\end{theorem}

\begin{proof}
The result follows directly from the variational principle applied to the functional 
$$
\int \Big(\sqrt{ |\langle \dot{\bf x}, \dot{\bf x}\rangle|}+\lambda( \langle {\bf x}, {\bf x}\rangle-1)\Big)dt ,
$$
where the dot stands for the derivative with respect to the parameter $t$ along the curve, with
corresponding Euler-Lagrange equations 
$$
{\ddot{\bf x}}=2\lambda {\bf x}, \quad \langle {\bf x}, {\bf x}\rangle=1 \label{geod}
$$
with $t$ is chosen so that $|\langle\dot{\bf x}, \dot{\bf x}\rangle|=1$.
Differentiating two times the second equation and combining with the first equation we find the expression
$$
\lambda=-\frac{\langle\dot{\bf x}, \dot{\bf x}\rangle}{2}.
$$
One shows that the $\lambda=-1/2$  solution can be ruled out, and in the case of $\lambda=1/2$, it is expressed in terms of hyperbolic functions
$$
{\bf x}={\bf u}\cosh s+{\bf v}\sinh s.
$$
Applying the conditions $-\langle\dot{\bf x}, \dot{\bf x}\rangle=\langle{\bf x}, {\bf x}\rangle=1$, we find that ${\bf u}, {\bf v}$ satisfy the conditions of the theorem. 
\end{proof}

\bigskip
For any ${\bf u}\in\hyperp$ and ${\bf v}\in\hyperl$ with $\langle {\bf u},{\bf v}\rangle =0$, we shall let 
$$L_{{\bf u},{\bf v}}=\{{\rm cosh}\, s~{\bf u}+{\rm sinh}\, s~{\bf v}:s\in{\mathbb R}\}$$ denote the supergeodesic determined in accordance with Theorem \ref{uvthm}.

\bigskip

\begin{remark} The relationship between the bodies ${\bf u}\in\hyperp', {\bf v}\in\hyperl'$ in the classical case and the usual dual ${\bf h}\in\hyperl'$ to the geodesic is as follows. The vectors ${\bf e}={\bf u}+{\bf v},{\bf f}={\bf u}-{\bf v}$ determine
${\bf d}\in{\mathbb R}^{2,1}$ so that $\langle{\bf d},{\bf d}\rangle=0$ and $\langle{\bf d},{\bf e}\rangle=\langle{\bf d},{\bf f}\rangle=2$.  This follows directly from considerations of the signature of the restriction of the inner product to the subspace spanned by ${\bf e},{\bf f}$, and indeed there are two such choices for ${\bf d}$ with ${\bf d},{\bf e},{\bf f}$ providing a basis for ${\mathbb R}^{3|0}$, one on either side of the hyperplane determined by ${\bf e}$ and ${\bf f}$.  One easily computes that ${\bf h}
=2^{-{1\over 2}}({\bf e}+{\bf f}-{\bf d})\in\hyperl'$ in this basis, so that $\langle{\bf h},{\bf u}\rangle=\langle{\bf h},{\bf v}\rangle=0$.
In effect, the geodesic which is the dual of ${\bf h}\in\hyperl'$ is the geodesic which is perpendicular to the dual of ${\bf v}\in\hyperl'$ that contains ${\bf u}\in\hyperp'$.  This ${\bf h}$ is especially significant in the classical setting since given two geodesics $L,L'$ with respective
${\bf h},{\bf h}'\in\hyperl'$, their square inner product $\langle {\bf h},{\bf h}'\rangle^2$ has geometric significance; see $\S$\ref{sec:2lines} for details and for contrast with the super case.
 \end{remark}

\begin{corollary} \label{geocor}
For any two distinct points ${\bf x}_1,{\bf x}_2\in\hyperp$, there is a unique supergeodesic containing them,
and the distance $d({\bf x}_1,{\bf x}_2)$ between them satisfies ${\rm cosh}\,d({\bf x}_1,{\bf x}_2)=\langle {\bf x}_1,{\bf x}_2\rangle$.
\end{corollary}
\begin{proof}
We begin with existence.  First note that $\langle{\bf x}_1,{\bf x}_2\rangle~>1$ for distinct
${\bf x}_1,{\bf x}_2\in\hyperp$ and define
$$\lambda=\sqrt{{\langle{\bf x}_1,{\bf x}_2\rangle+1}\over{\langle{\bf x}_1,{\bf x}_2\rangle-1}}.$$
It follows that $\langle{\bf x}_1,{\bf x}_2\rangle={{\lambda^2+1}\over{\lambda^2-1}}$, whence
${\rm cosh}^{-1}t={\rm ln}(t+\sqrt{t^2-1})$ gives 
$$e^{{\rm cosh}^{-1}\langle{\bf x}_1,{\bf x}_2\rangle}=\langle{\bf x}_1,{\bf x}_2\rangle+\sqrt{\langle{\bf x}_1,{\bf x}_2\rangle^2-1}={{\lambda+1}\over{\lambda -1}}.$$

In the notation of
Theorem \ref{uvthm}, we directly exhibit $L_{{\bf u},{\bf v}}=L_{{{{\bf e}+{\bf f}}\over 2},{{{\bf e}-{\bf f}}\over 2}}$ containing ${\bf x}_1,{\bf x}_2$, where
$${\bf x}(s)={\rm cosh}\, s~{\bf u}+{\rm sinh}\, s~{\bf v}={1\over 2}(e^s~{\bf e}+e^{-s}~{\bf f})$$
taking
$$\begin{aligned}
{\bf e}&={{\lambda-1}\over{2\lambda}}[(1-\lambda){\bf x}_1+(1+\lambda){\bf x}_2],\\
{\bf f}&={{\lambda+1}\over{2\lambda}}[(1+\lambda){\bf x}_1+(1-\lambda){\bf x}_2].\\
\end{aligned}$$
It is easy to confirm directly that ${\bf e},{\bf f}\in{\mathcal L}^+$ with $\langle{\bf e},{\bf f}\rangle=2$ and that
$${\bf x}(0)={\bf x}_1~{\rm and}~{\bf x}({\rm cosh}^{-1}\langle{\bf x}_1,{\bf x}_2\rangle)={\bf x}_2,$$
thus establishing existence as well as that ${\rm cosh}\,d({\bf x}_1,{\bf x}_2)=\langle {\bf x}_1,{\bf x}_2\rangle$. 

\bigskip

As for uniqueness, the system of equations
$$\begin{aligned}
{\bf x}_1&={\rm cosh}\,p~{\bf u}+{\rm sinh}\, p~{\bf v},\\
{\bf x}_2&={\rm cosh}\,q~{\bf u}+{\rm sinh}\, q~{\bf v}
\end{aligned}$$
is tantamount to the linear system
$$\begin{pmatrix}{\rm cosh}\, p\,I&{\rm sinh}\, p\,I\\{\rm cosh}\, q\,I&{\rm sinh}\, q\,I\end{pmatrix}~\begin{pmatrix}{\bf u}\\{\bf v}\end{pmatrix}=\begin{pmatrix}{\bf x}_1\\{\bf x}_2\end{pmatrix},$$
where $I$ is the 5-by-5 identity matrix, which has a unique solution since the determinant is non-zero.
for $p\neq q$.\end{proof}

\begin{corollary}\label{useful}
Let $L=L_{{\bf u},{\bf v}}$ and set ${\bf e}={\bf u}+{\bf v}, {\bf f}={\bf u}-{\bf v}\in{\mathcal L}^+$.
Then we have
$$\begin{aligned}
L&=\{ \bf P\in\hyperl:\langle {\bf P},{\bf e}\rangle~\langle {\bf P},{\bf f}\rangle~=~1\}\, ,\\
&=\biggl\{ {1\over{2\sqrt{xy}}}(x{\bf e}+y{\bf f}):x,y> 0\biggr\}\,.\\
\end{aligned}$$
\end{corollary}

\begin{proof} For the inclusion of $L$ in the first equality, write $${\bf P}={\rm cosh}\, p~{\bf u}~+~{\rm sinh}\, p~{\bf v}\,,$$ so that $\langle {\bf P},{\bf u}\pm{\bf v}\rangle={\rm cosh}\, p\mp{\rm sinh}\,p,$ so
$$\langle {\bf P},{\bf e}\rangle~\langle {\bf P},{\bf f}\rangle={\rm cosh}^2 p-{\rm sinh}^2 p=1.$$

\bigskip

For the reverse inclusion, suppose that
$$\begin{aligned}
1&=\langle {\bf P},{\bf e}\rangle~\langle {\bf P},{\bf f}\rangle\\
&=\langle {\bf P},{\bf u}\rangle^2-\langle {\bf P},{\bf v}\rangle^2\, \\
\end{aligned}$$
and define ${\bf Q}=\langle {\bf P},{\bf u}\rangle\,{\bf u}-\langle {\bf P},{\bf v}\rangle\,{\bf v}$,
whence
$$\begin{aligned}
\langle{\bf Q},{\bf Q}\rangle&=\biggl\langle \langle {\bf P},{\bf u}\rangle\,{\bf u}-\langle {\bf P},{\bf v}\rangle\,{\bf v}~,~\langle {\bf P},{\bf u}\rangle\,{\bf u}-\langle {\bf P},{\bf v}\rangle\,{\bf v}\biggr\rangle\\
&=\langle{\bf P},{\bf u}\rangle^2-\langle{\bf P},{\bf v}\rangle^2\\
&=1,
\end{aligned}$$
so ${\bf Q}\in\hyperp$.  Moreover by the previous corollary, the cosh of the distance between
${\bf P}$ and ${\bf Q}$ is given by
$$\begin{aligned}
\langle{\bf P},{\bf Q}\rangle&=\biggr\langle{\bf P}, \langle {\bf P},{\bf u}\rangle\,{\bf u}-\langle {\bf P},{\bf v}\rangle\,{\bf v}\biggr\rangle\\
&=\langle{\bf P},{\bf u}\rangle^2-\langle{\bf P},{\bf v}\rangle^2,\\
\end{aligned}$$
so the distance between ${\bf P}$ and ${\bf Q}$ is zero, whence ${\bf P}={\bf Q}$,
proving the first identity.

\bigskip

For the second equality, suppose
${\bf Q}=x{\bf e}+y{\bf f}$, for $x,y>0,$ so
$$\langle {\bf Q},{\bf e}\rangle~\langle{\bf Q},{\bf f}\rangle~=~4xy=\langle{\bf Q},{\bf Q}\rangle$$
since $\langle{\bf e},{\bf f}\rangle=2$.  The first equality shows that 
${1\over{2\sqrt{xy}}}{\bf Q}\in L_{{\bf u},{\bf v}}$ and the second that ${1\over{2\sqrt{xy}}}{\bf Q}\in\hyperp$ as required.
\end{proof}

\noindent The significant computational importance of this result is that a supergeodesic is the projectivization to lie in $\hyperp$ of the convex linear span of vectors lying in its
asymptotes in ${\mathcal L}^+$.

\bigskip

\begin{proposition}\label{prop:luvs} Let ${\bf e}_i,{\bf f}_i\in{\mathcal L}^+$ with
$\langle {\bf e}_i,{\bf f}_i\rangle=2$ and define
$${\bf u}_i={{{\bf e}_i+{\bf f}_i}\over 2},~~{\bf v}_i={{{\bf e}_i-{\bf f}_i}\over 2}, ~~{\rm for}~i=1,2.$$ Then $L_{{\bf u}_1,{\bf v}_1}=L_{{\bf u}_2,{\bf v}_2}$ as oriented supergeodesics if and only if
${\bf e}_1$,${\bf e}_2$ (and ${\bf f}_1$,${\bf f}_2$, respectively) are proportional.
Moreover in this case ${\bf e}_2={\rm ln}\,a~{\bf e}_1$, for $a\geq 1$, implies
${\bf f}_2={1\over {{\rm ln}\, a}}{\bf f}_1$, and 
$${\rm cosh}\, p~{\bf u}_1~+~{\rm sinh}\, p~{\bf v}_1=
{\rm cosh}(p-a)~{\bf u}_2~+~{\rm sinh}(p-a)~{\bf v}_2,~~{\rm for~all}~s.$$
\end{proposition}

\bigskip

In other words, scaling ${\bf u}+{\bf v}$ and ${\bf u}-{\bf v}$ by reciprocal amounts merely shifts the origin of the parametrization of $L_{{\bf u},{\bf v}}$, and interchanging them reverses the orientation.  We shall refer to the point 
${\bf u}\in L_{{\bf u},{\bf v}}$ arising from vanishing parameter as the {\it origin} of the parametrization, so ${\bf u},{\bf v}$ determines not only the line but also an origin within it.

\begin{proof}
In case $L_{{\bf u}_1,{\bf v}_1}=L_{{\bf u}_2,{\bf v}_2}$, the claimed proportionality follows easily from the second part of the previous result.  For the converse with
$b={\rm ln}\, a$, ${\bf e}_2=b{\bf e}_1$ and  ${\bf f}_2=c{\bf f}_1$ implies
$$2=\langle {\bf e}_2,{\bf f}_2\rangle=\langle b{\bf e}_1,c{\bf f}_1\rangle=2bc,$$
whence $bc=1$.  Thus, we find
$$\begin{aligned}
{\bf u}_2&=\mbox{\normalsize${{b+b^{-1}}\over 2}$}{\bf u}_1+\mbox{\normalsize${{b-b^{-1}}\over 2}$}{\bf v}_1={\rm cosh}\, a~{\bf u}_1+{\rm sinh}\, a~{\bf v}_1,\\
{\bf v}_2&=\mbox{\normalsize${{b-b^{-1}}\over 2}$}{\bf u}_1+\mbox{\normalsize${{b+b^{-1}}\over 2}$}{\bf v}_1={\rm sinh}\, a~{\bf u}_1+{\rm cosh}\, a~{\bf v}_1,\\
\end{aligned}$$
and so
$$\begin{aligned}
{\rm cosh}\, p~{\bf u}_1+{\rm sinh}\, p~{\bf v}_1&=
{\rm cosh}\, q\,[{\rm cosh}\, a~{\bf u}_1+{\rm sinh}\, a~{\bf v}_1]\\
&\hskip.6ex+{\rm sinh}\, q\,[{\rm sinh}\, a~{\bf u}_1+{\rm cosh}\, a~{\bf v}_1]\\
\end{aligned}$$
follows from
$$\begin{aligned}
{\rm cosh}\, p&={\rm cosh}\, a~{\rm cosh}\, q~+~{\rm sinh}\, a~{\rm sinh}\, q={\rm cosh}\,(a+q),\\
{\rm sinh}\, p&={\rm sinh}\, a~{\rm cosh}\, q~+~{\rm cosh}\, a~{\rm sinh}\, q={\rm sinh}\,(a+q),\\
\end{aligned}$$
upon taking respective inner products with ${\bf u}_1$, ${\bf v}_1$, whence
$${\rm sinh}\, 2p={{{\rm cosh}\,p~{\rm sinh}\,p}\over 2}=
{{{\rm cosh}\,(a+q)~{\rm sinh}\,(a+q)}\over 2}={\rm sinh}\, 2(a+q),$$
and so $p =a+q$.
\end{proof}

\section{Law of Cosines}

It follows from Corollary \ref{geocor} that any three points of $\hyperp$ not lying on a common supergeodesic determine a
{\it supertriangle}, namely, three geodesic segments with disjoint interiors meeting pairwise
at the given points.  We adopt the usual terminology of sides and opposite angles from their bodies, and have

\begin{corollary}\label{LOC}
The usual hyperbolic laws of cosines and sines hold for supertriangles, that is
given a supertriangle with sides of edge lengths $A,B,C$ and opposite angles $a,b,c$, then we have
$$\begin{aligned}
{\rm cosh}\,A&={\rm cosh}\,B ~{\rm cosh}\,C~-~{\rm sinh} \,B ~{\rm sinh} \,C~ {\rm cos\,}a,~~~%{\rm (HLoCI),}
\\
{\rm cos}\,a&=-{\rm cos}\,b~{\rm cos}\,c~+~{\rm sin}\,b~{\rm sin}\,c~{\rm cosh}\,A,
%\hskip .7cm{\rm (HLoCII),}
\\
\frac{{\rm sin}\,a}{{\rm sinh}\,A}&=\frac{{\rm sin}\,b}{{\rm sinh}\,B}=\frac{{\rm sin}\,c}{{\rm sinh}\,C}.%\hskip 4.0cm{\rm (HLoS)}.
\\
\end{aligned}$$
\end{corollary}
\begin{proof}
According to Theorem \ref{uvthm}, the supergeodesic from ${\bf X}\in\hyperp$ to ${\bf Y}\in\hyperp$ may be parametrized
$${\bf X}~{\rm cosh}\,t+{{{\bf Y}-{\bf X}\,\langle {\bf X},{\bf Y}\rangle}\over{\sqrt{\langle {\bf X},{\bf Y}\rangle^2-1}}}\,{\rm sinh}\,t,~~{\rm for}~ t\in{\mathbb R},$$
 where the unit tangent vector at ${\bf X}$ to the line from ${\bf X}$ to ${\bf Y}$ is given by
 ${{\bf Y}-{\bf X}\,\langle {\bf X},{\bf Y}\rangle}\over{\sqrt{\langle {\bf X},{\bf Y}\rangle^2-1}}$.  Letting ${\bf X}$ denote the vertex opposite the side of length $X$, for $X=A,B,C$ and taking the inner
 product of the unit tangent vectors at ${\bf A}$ to the supergeodesics through ${\bf A}$ and ${\bf B},{\bf C}$ thus gives
$$
{\rm cos}\,a=\biggr\langle {{{\bf B}-{\bf A}\,\langle {\bf A},{\bf B}\rangle}\over{\sqrt{\langle {\bf A},{\bf B}\rangle^2-1}}},
{{{\bf C}-{\bf A}\,\langle {\bf A},{\bf C}\rangle}\over{\sqrt{\langle {\bf A},{\bf C}\rangle^2-1}}}\biggr\rangle,
$$
and so 
$$\begin{aligned}
{\rm cos}\,a~{\rm sinh}\,B~{\rm sinh}\, C
&=\langle {\bf B},{\bf C}\rangle -\langle {\bf A},{\bf B}\rangle {\rm cosh}\,B-\langle {\bf A},{\bf C}\rangle {\rm cosh}\,C\\
&~~~+\langle{\bf A},{\bf A}\rangle~{\rm cosh}\,B~{\rm cosh}\, C\\
&={\rm cosh}\, A-{\rm cosh}\,C~{\rm cosh}\, B-{\rm cosh}\, B~{\rm cosh}\, C\\
&~~~+{\rm cosh}\,B~{\rm cosh}\,C\\
&={\rm cosh}\,A-{\rm cosh}\, B~{\rm cosh}\, C,
\end{aligned}$$
proving the first identity.  It is a well-known and easy exercise to derive the second and third identities purely algebraically from this using just ${\rm cos}^2x+{\rm sin}^2x=1={\rm cosh}^2X-{\rm sinh}^2X$.
\end{proof}

\begin{corollary} Given a supergeodesic $L$ and a point ${\bf P}\notin L$,
there exists a unique supergeodesic $L'$ through ${\bf P}$ perpendicular to $L$.
\end{corollary}

\begin{proof} Suppose that $L=L_{{\bf u},{\bf v}}$, and let $d=d(p)>0$ denote the distance from ${\bf P}$ to the
point ${\rm cosh}\, p~{\bf u}~+~{\rm sinh}\, p~{\bf v}\in L$.  According to  Corollary \ref{geocor}, we have
$$\begin{aligned}
{\rm cosh}\, d&=\langle {\bf P},{\rm cosh}\, p~{\bf u}~+~{\rm sinh}\, p~{\bf v}\rangle\\
&={\rm cosh}\,p ~\langle {\bf P},{\bf u}\rangle+{\rm sinh}\,p ~\langle {\bf P},{\bf u}\rangle\,.
\end{aligned}$$
The unique critical point ${{d}\over{dp}}{\rm cosh}\, d=0$ occurs for 
$${\rm tanh}\, p=-{{\langle {\bf P},{\bf u}\rangle}\over{\langle {\bf P},{\bf v}\rangle}},$$
which is a minimum since ${{d^2}\over{dp^2}}{\rm cosh}\, d={\rm cosh}\, d>0.$
Uniqueness follows from convexity of {\rm cosh} for positive argument, and the Law of Cosines easily implies that this minimizer is perpendicular to $L$.
\end{proof}

\noindent The point 
$$\begin{aligned}
{\rm cosh}\,p~{\bf u}+{\rm sinh}\, p~{\bf v}
&={{ {\bf u}+p~{\bf v}}\over\sqrt{p^2+1}},~{\rm for}~{\rm tanh}\, p=-{{\langle {\bf P},{\bf u}\rangle}\over{\langle {\bf P},{\bf v}\rangle}}\\
&={{\langle {\bf P},{\bf v}\rangle\,{\bf u}~-~\langle {\bf P},{\bf u}\rangle\,{\bf v}}\over
\sqrt{\langle {\bf P},{\bf u}\rangle^2+\langle {\bf P},{\bf v}\rangle^2}}\\
\end{aligned}$$
is the {\it orthogonal projection} of ${\bf P}$ on $L_{{\bf u},{\bf v}}$.

\bigskip

\section{Pairs of supergeodesics}\label{sec:2lines}

\noindent We begin with a technical lemma.

\begin{lemma}
Assume that
$$a_1{\bf d}+b_1{\bf e}+c_1{\bf f}+(0,0,0~\vert~\alpha_1,\beta_1)=
a_2{\bf d}+b_2{\bf e}+c_2{\bf f}+(0,0,0~\vert~\alpha_2,\beta_2),$$
where ${\bf d}=(*,*,*~\vert~0,0),~{\bf e}=(*,*,*~\vert~\phi,0),~
{\bf f}=(*,*,*~\vert~0,\psi)\in{\mathcal L}$
with $\langle {\bf d},{\bf e}\rangle=\langle {\bf e},{\bf f}\rangle=
\langle {\bf f},{\bf d}\rangle=2$. Then
$a_1=a_2, b_1=b_2, c_1=c_2$ and $\alpha_1=\alpha_2,\beta_1=\beta_2$.
\end{lemma}
\begin{proof}
The respective inner products of the assumed equality with ${\bf d},{\bf e}$, ${\bf f}$ yield
$M\begin{psmallmatrix}a_2\\b_2\\c_2\end{psmallmatrix}=M\begin{psmallmatrix}a_2\\b_2\\c_2\end{psmallmatrix}+\begin{psmallmatrix}0\\\phi(\beta_1-\beta_2)\\(\alpha_1-\alpha_2)\psi\\\end{psmallmatrix}$, where $M=\begin{psmallmatrix}0&1&1\\1&0&1\\1&1&0\\\end{psmallmatrix},$ so that
$$\begin{pmatrix}a_2\\b_2\\c_2\\\end{pmatrix}=\begin{pmatrix}a_1\\b_1\\c_1\\\end{pmatrix}
+{1\over 2}\begin{pmatrix}-1&\hskip 1.8ex1&\hskip 1.8ex1\\\hskip 1.8ex1&-1&\hskip 1.8ex1\\\hskip 1.8ex1&\hskip 1.8ex1&-1\\\end{pmatrix}
\begin{pmatrix}0\\\phi(\beta_1-\beta_2)\\(\alpha_1-\alpha_2)\psi\\\end{pmatrix}$$
and with $(0,0,0~\vert~\phi,0)$ and $(0,0,0~\vert~0,\psi)$ yield
$$\begin{aligned}
(\alpha_1-\alpha_2)\psi&=(b_2-b_1)\phi\psi,\\
\phi(\beta_1-\beta_2)&=(c_2-c_1)\phi\psi.\\
\end{aligned}$$
It follows that
$$\begin{aligned}
a_2&=a_1+{1\over 2}[(b_2-b_1)+(c_2-c_1)]\phi\psi,\\
b_2-b_1&={1\over 2}[(b_2-b_1)-(c_2-c_1)]\phi\psi,\\
c_2-c_1&={1\over 2}[(c_2-c_1)-(b_2-b_1)]\phi\psi,\\
\end{aligned}$$
whence $a_1=a_2$, and likewise for $b_1=b_2$, $c_1=c_2$, from which it
finally follows that $\alpha_1=\alpha_2$, $\beta_1=\beta_2$.
\end{proof}

\begin{theorem}\label{thminter}
Consider the supergeodesic $L=L_{{\bf u},{\bf v}}$ determined by
$$\begin{aligned}
{\bf u}&=(*,*,*~\vert~{\phi\over 2}, +{\psi\over 2})\in\hyperp,\\
{\bf v}&=(*,*,*~\vert~{\phi\over 2},-{\psi\over 2})\in\hyperl\,.
\end{aligned}$$
Let 
$$
{\bf e}={\bf u}+{\bf v}=(*,*,*~\vert~\phi,0),~~
{\bf f}={\bf u}-{\bf v}=(*,*,*~\vert~0,\psi)\in{\mathcal L}^+,$$
and choose $${\bf d}=(*,*,*~\vert~0,0)\in{\mathcal L}^+,~{\rm so~that}
~\langle{\bf d},{\bf e}\rangle = \langle{\bf e},{\bf f}\rangle=\langle{\bf f},{\bf e}\rangle=2.$$
Suppose that $L'=L_{{\bf u}',{\bf v}'}$ is another supergeodesic and write 
$$\begin{aligned}
{\bf u}'&=a{\bf d}+b{\bf e}+c{\bf f}+(0,0,0~\vert~2\alpha,2\beta),\\
{\bf v}'&=x{\bf d}+y{\bf e}+z{\bf f}+(0,0,0~\vert~2\xi,2\eta).\\
\end{aligned}$$
Define
$$\begin{aligned}
A&=\langle {\bf u},{\bf u}'\rangle,\quad B=\langle {\bf v},{\bf u}'\rangle,\\
C&=\langle {\bf u},{\bf v}'\rangle,\quad D=\langle {\bf v},{\bf v}'\rangle,\\
\end{aligned}$$
and set 
$$\begin{aligned}
I&=[A^2-B^2-1]^{1\over 2}~\biggl [1+{{\alpha\beta(4+\phi\psi)}\over{A^2-B^2-1}}\biggr ],\\
J&=[C^2-D^2+1]^{1\over 2}~\biggl [1+{{\xi\eta(4+\phi\psi)}\over{C^2-D^2+1}}\biggr ].\\
\end{aligned}$$
Then $L$ and $L'$ intersect if and only if the following conditions hold:
\medskip

\leftskip 1.5cm

\noindent
{\rm (\ref{thminter}.1)}~$C^2-D^2+1\geq A^2-B^2-1\geq 0$,\\

\noindent 
{\rm (\ref{thminter}.2)}~$AC-BD+IJ=(\eta\alpha+\beta\xi)(4+\phi\psi)$,\\

\noindent {\rm (\ref{thminter}.3)} $-2IJ(AC-BD)=J^2(A^2-B^2-1)+I^2(C^2-D^2+1)$,\\

\noindent 
{\rm (\ref{thminter}.4)}~$4(J\alpha+I\xi)=[(A-B)J+(C-D)I]\phi$,\\

\noindent 
{\rm (\ref{thminter}.5)}~$4(J\beta+I\eta)=[(A+B)J+(C+D)I]\psi$.\\

\leftskip=0ex

Moreover, under these conditions, the unique point of intersection is given by
$$
L\cap L'={{J{\bf u}'+I{\bf v}'}\over[J^2-I^2]^{{1\over 2}}} 
={{(AJ+CI){\bf u}+(BJ+DI){\bf v}}\over[J^2-I^2]^{{1\over 2}}}.$$
\end{theorem}

\begin{proof}
It follows from the definitions that
$$\begin{aligned}
A&=2a+b+c+\alpha\psi+\phi\beta,~~B=c-b-\alpha\psi+\phi\beta,\\
C&=2x+y+z+\xi\psi+\phi\eta,~~D=z-y-\xi\psi+\phi\eta,\\
\end{aligned}$$
whence
$$\begin{aligned}
b=-a+{{A-B}\over 2}-\alpha\psi,~~c=-a+{{A+B}\over 2}-\phi\beta,\\
y=-x+{{C-D}\over 2}-\xi\psi,~~z=-x+{{C+D}\over 2}-\phi\eta.\\
\end{aligned}$$
Moreover, we have
\begin{eqnarray}
1&=&\langle{\bf u}',{\bf u}'\rangle=4(ab+bc+ca)+4(b\phi\beta+c\alpha\psi+2\alpha\beta),\label{eqq1}\\
-1&=&\langle{\bf v}',{\bf v}'\rangle=4(xy+yz+zx)+4(y\phi\eta+z\xi\psi+2\xi\eta),\label{eqq2}\\
0&=&\langle{\bf u}',{\bf v}'\rangle=2[a(y+z)+b(x+z)+c(x+y)]\label{eqq3}\\
&~+&2(c\xi\psi+b\phi\eta+z\alpha\psi+y\phi\beta)+4(\alpha\eta+\xi\beta).\nonumber
\end{eqnarray}
Now writing eqn (\ref{eqq1}) in terms of $a$ alone gives the quadratic
$$0=-4a^2-4a(\alpha\psi+\phi\beta)+[A^2-B^2-1+4\alpha\beta(2+\phi\psi)],$$
from which it follows that
$$a=-{1\over 2}(\alpha\psi+\phi\beta\pm I),$$
where $I$ is defined in the statement of the theorem.

\medskip

The analogous computation using eqn (\ref{eqq2}) gives 
$$x=-{1\over 2}(\xi\psi+\phi\eta\pm J),$$
yielding first of all the two constraints $A^2-B^2-1\geq 0\leq C^2-D^2+1$ which form part of condition (\ref{thminter}.1) in the theorem.  Moreover, $a$ and $x$ must have opposite signs in order for $L$ to intersect $L'$, and the sign of $I$ in $a$ must be positive in order that ${\bf u}'\in\hyperp$.  We therefore find
$$a=-{1\over 2}(\alpha\psi+\phi\beta- I)~{\rm and}~x=-{1\over 2}(\xi\psi+\phi\eta + J).$$

\medskip

Writing eqn  (\ref{eqq3}) in terms of $a$ and $x$ gives
$$\begin{aligned}
AC-B&D+4(\alpha\eta+\xi\beta)+2\phi\psi(\alpha\eta+\xi\beta)\\=&~~4ax+2a(\phi\eta+\xi\psi)+2x(\alpha\psi+\phi\beta),\\
\end{aligned}$$
and finally plugging in the values of $a$ and $x$ yields condition
(\ref{thminter}.2) of the theorem.

\medskip

Now take the respective inner products of 
\begin{eqnarray}\rm cosh}\,p~{\bf u}'
+~{\rm sinh}\,p~{\bf v}'~~=~~{\rm cosh}\,q~{\bf u}~+~{\rm sinh}\,q~{\bf v}{\label{equals}
\end{eqnarray}
with ${\bf d},{\bf e},{\bf f}$ to get the linear system
$$MX\begin{pmatrix}a\\b\\c\\\end{pmatrix}+MY\begin{pmatrix}x\\y\\z\\\end{pmatrix}
~~=~~\begin{pmatrix}2{\rm cosh}\,q\\{\rm cosh}\,q-{\rm sinh}\, q\\{\rm cosh}\,q+{\rm sinh}\, q\\\end{pmatrix}-2\begin{pmatrix}0\\\phi\beta~{\rm cosh}\,p+\phi\eta~{\rm sinh}\,p\\
\alpha\psi~{\rm cosh}\,p+\xi\psi~{\rm sinh}\,p\end{pmatrix}\\,$$
where $X$ and $Y$ are 3-by-3 diagonal matrices with respective diagonal entries 
${\rm cosh}\, s$ and ${\rm sinh}\, s$ with 
$$M=\begin{pmatrix}0&2&2\\2&0&2\\2&2&0\\\end{pmatrix} ~{\rm so~that}~
M^{-1}={1\over 4}\begin{pmatrix}-1&\hskip 1.9ex1&\hskip 1.9ex1\\\hskip 1.9ex1&-1&\hskip 1.9ex1\\\hskip 1.9ex1&\hskip 1.9ex1&-1\\\end{pmatrix}.$$
It follows that
\begin{eqnarray}
\label{eqni}2(a{\rm cosh}\,p+x{\rm sinh}\, p)&=&(\alpha\psi+\phi\beta){\rm cosh}\,p+(\xi\psi+\phi\eta){\rm sinh}\, p ,~~~\\
\label{eqnii}2(b{\rm cosh}\,p+y{\rm sinh}\, p)&=&(\phi\beta-\alpha\psi){\rm cosh}\,p+(\phi\eta-\xi\psi){\rm sinh}\, p \nonumber\\
~~&+~~~&{\rm cosh}\, q+{\rm sinh}\, q,\\
\label{eqniii}2(c{\rm cosh}\,p+z{\rm sinh}\, p)&=&(\alpha\psi-\phi\beta){\rm cosh}\,p+(\xi\psi-\phi\eta){\rm sinh}\, p \nonumber\\
~~&~~~+&{\rm cosh}\, q-{\rm sinh}\, q.
\end{eqnarray}
Eqn (\ref{eqni}) gives
$${\rm tanh}\,p={{{\rm sinh}\, p}\over{{\rm cosh}\, p}}={I\over J},$$
whence
$$
{\rm cosh}\, p={J\over{[J^2-I^2]^{1\over 2}}}~{\rm and}~
{\rm sinh}\, p={I\over{[J^2-I^2]^{1\over 2}}},
$$
taking the parameter $p\geq 0$.  Eqn (\ref{eqni}) therefore furthermore implies the inequality $C^2-D^2+1\geq A^2-B^2-1$,
thus completing the proof of necessity of condition (\ref{thminter}.1) of the theorem.

\medskip

Meanwhile, the sum of eqn (\ref{eqnii}) and eqn (\ref{eqniii}) provides
$$\begin{aligned}
{\rm cosh}\, q&=(b+c){\rm cosh}\, p+(y+z){\rm sinh}\, p\\
&=(A-I){\rm cosh}\,p +(C+J){\rm sinh}\, p,\\
\end{aligned}$$
while their difference yields
$$\begin{aligned}
{\rm sinh}\, q&=[b-c+(\alpha\psi-\phi\beta)]{\rm cosh}\, p+[y-z+(\xi\psi-\phi\eta)]{\rm sinh}\,p\\
&=-B{\rm cosh}\, p-D{\rm sinh}\, p.\\
\end{aligned}$$
Computation shows that the constraint $1={\rm cosh}^2\, q-{\rm sinh}^2\, q$ is equivalent to
condition (\ref{thminter}.3) of the theorem.

\medskip

The final two fermionic constraints (\ref{thminter}.4) and (\ref{thminter}.5) arise by equating the last two coordinate entries of eqn (\ref{equals}).  The stated necessary conditions are clearly sufficient, and the formula for $L\cap L'$ then follows immediately from eqn (\ref{equals}).
\end{proof}

\begin{remark}
For any ${\bf e},{\bf f}\in{\mathcal L}^+$ with $\langle {\bf e},{\bf f}\rangle=2$, it is an easy matter to find $g\in{\rm OSp}(1\vert 2)$ so that $g.{\bf e}$ and $g.{\bf f}$ satisfy the conditions of the previous theorem.  Indeed for ${\bf e}=(x_1,x_2,y~\vert~\xi,\eta)$, ${\bf f}=(u_1,u_2,v~\vert~\mu,\nu)$,
we may simply take $g=u\bigl({{x_1\nu-v\xi}\over{yv-x_1u_2}},
{{y\nu-u_2\xi}\over{yv-x_1u_2}}\bigr)$ in the generic case that $yv\neq x_1u_2$, which can itself also be easily arranged by explicit perturbation.
Since $g(L\cap L')=g(L)\cap g(L')$, for any two supergeodesics $L,L'$ and $g\in{\rm OSp}(1\vert2)$, the previous result in fact gives conditions and formulae for intersections in the general case.
\end{remark}

\begin{corollary}\label{cosines}
Consider supergeodesics $L=L_{{\bf u},{\bf v}}$ and $L'=L_{{\bf u}',{\bf v}'}$ in the notation of Theorem \ref{uvthm} and set
$X=A^2-B^2-1,~~Y=C^2-D^2+1$
with
$Z=Y+X+2(AB+CD)(BD-AC)$.
If $L$ and $L'$ intersect at the point $P\in\hyperp$, then the cosine of the angle at $P$ from $L$ to $L'$ is given by
$${Z\over{Y-X}}+{{2(4+\phi\psi)}\over{Y-X}}\biggl [
{{(\alpha\beta+\xi\eta)+(AB+CD)(\eta\alpha+\beta\xi)}\atop
{+{{Z(\alpha\beta-\xi\eta)}\over{Y-X}}+{{16Z\alpha\beta\xi\eta}\over{(Y-X)^2}}}}
\biggr ].$$
\end{corollary}

\medskip

\noindent More interesting than the particular form of the second summand is 
the substantial fermionic correction to the bosonic term ${Z\over{Y-X}}$ that it
represents.  This puts into focus the Super Law of Cosines Corollary \ref{LOC}, where
the familiar purely bosonic formula rather remarkably concisely includes this and
other fermionic corrections.

\medskip

\begin{proof}
From the formulae for $P=L\cap L'$ in Theorem \ref{uvthm}, the cosine of the angle at $P$ from $L$ to $L'$ is given by
$$\begin{aligned}
&{{\langle J{\bf v}'+I{\bf u}',(AJ+CI){\bf v}+(BJ+DI){\bf u}\rangle}\over{J^2-I^2}}\\
&\hskip.5cm={{(AD+BC)(J^2+I^2)+2(AB+CD)IJ}\over{J^2-I^2}}.
\end{aligned}$$
Plugging in the constraint 
$IJ=BD-AC+(\eta\alpha+\beta\xi)(4+\phi\psi)$
from {\rm (\ref{thminter}.2)} gives the asserted formula after a bit of computation
using constraints {\rm (\ref{thminter}.3)} and {\rm (\ref{thminter}.4)} to conclude
that $\alpha\beta\xi\eta\phi\psi=0$.
\end{proof}

A pair of supergeodesics $L_i=L_{{{\bf u}_i},{{\bf v}_i}}$, for $i=1,2$, in the bosonic case are determined by a pair ${\bf u}_i\in\hyperp'$ and ${\bf v}_i\in\hyperl'$ with $\langle {\bf u}_i,{\bf v}_i\rangle=0$,
or equivalently by a pair ${\bf h}_i\in\hyperl$, where $\langle {\bf h}_i,{\bf u}_i\rangle=\langle {\bf h}_i,{\bf v}_i\rangle=0$, and in this case \cite{Thurston} we have

$$\langle {\bf h_1},{\bf h_2}\rangle^2 =\begin{cases}
{\rm cosh}^2\, d,&{\rm if}~L_1\cap L_2=\emptyset~{\rm and~are~not~ultraparallel},\\
&{\rm where}~d~{\rm is~the~distance~between}~L_1~{\rm and}~L_2;\\
{\rm cos}^2\, \alpha,&{\rm if}~L_1\cap L_2\neq\emptyset,\\
&{\rm where}~\alpha~{\rm is~the~angle~between}~L_1~{\rm and}~L_2,\\
%0,&{\rm if}~L_1~{\rm and}~L_2~{\rm are~ultraparallel}.
\end{cases}$$

\noindent where two geodesics or two supergeodesics
are {\it parallel} provided they do not intersect
and are furthermore {\it ultraparallel} if they share an asymptotic ray in ${\mathcal L}^+$.

\medskip

\noindent 
Moreover in the first case that $L_1$ and $L_2$ are parallel and not ultraparallel, then they admit a unique common perpendicular.  The sequel is dedicated to the analogue in hyperbolic superspace.

\begin{theorem} Two supegeodesics that are parallel and not
ultraparallel admit a unique common perpendicular. 
\end{theorem}

\begin{proof}
Denote the two supergeodesics $L_i=L_{{{\bf u}_i},{{\bf v}_i}}$, for $i=1,2$, and
define
$$\begin{aligned}
a&=\langle {\bf u}_1,{\bf u_2}\rangle,~~b=\langle {\bf v}_1,{\bf u}_2\rangle,\\
c&=\langle {\bf u}_1,{\bf v_2}\rangle,~~b=\langle {\bf v}_1,{\bf v}_2\rangle.\\
\end{aligned}$$
First note that the condition $\{{\bf u}_1\pm{\bf v_1}\}\cap\{{\bf u}_2\pm{\bf v_2}\}=\emptyset$ that $L_1$ and $L_2$ are not ultraparallel
thus implies positivity of each of the four linear combinations
$$(a-b-c+d), (a-b+c-d),(a+b-c-d),(a+b+c+d),$$
since the inner product of any two points in ${\mathcal L}^+$ which lie in different rays from the origin in ${\mathbb R}^{2,1|2}$ is positive.
In particular, the various pairwise sums of these inequalities imply that $a>{\rm max}\{ \pm b, \pm c,\pm d\}$, whence 
$a^2>{\rm max}\{b^2,c^2,d^2\}$.

\bigskip

Distinct points ${\bf u}_i{\rm cosh}\, t_i + {\bf v}_i{\rm sinh}\, t_i \in L_i$, for $i=1,2$, determine a unique supergeodesic $L$
containing them according to  Corollary \ref{geocor}, and furthermore the cosh of the length $d$ of the segment between $L\cap L_1$ and $L\cap L_2$ is given
by
$$\begin{aligned}
{\rm cosh}\, d&=\langle {\bf u}_1{\rm cosh}\, p_1 + {\bf v}_1{\rm sinh}\, p_1,{\bf u}_2{\rm cosh}\, p_2 + {\bf v}_2{\rm sinh}\, p_2\rangle\\
&=a~{\rm cosh}\,p_1~{\rm cosh}\, p_2~+~b~{\rm sinh}\,p_1~{\rm cosh}\, p_2~\\
&+~c~{\rm cosh}\,p_1~{\rm sinh} \,p_2~+~d~{\rm sinh}\,p_1~{\rm sinh}\, p_2,\\
\end{aligned}$$
whence
$$\begin{aligned}
{{\partial{\rm cosh}\,d}\over{\partial p_1}}&={\rm sinh}\, p_1\{a~{\rm cosh}\,p_2~+~c~{\rm sinh}\,p_2\} \\
&+~{\rm cosh}\, p_1\{b~{\rm cosh}\,p_2~+~d~{\rm sinh}\,p_2\} ,\\\\
{{\partial{\rm cosh}\,d}\over{\partial p_2}}&={\rm sinh}\, p_2\{a~{\rm cosh}\,p_1~+~b~{\rm sinh}\,p_1\} \\
&+~{\rm cosh}\, p_2\{c~{\rm cosh}\,p_1~+~d~{\rm sinh}\,p_1\} ,\\
\end{aligned}$$
and it follows that at a critical point
\begin{eqnarray}
\nonumber-{\rm tanh}\, p_1= {{d~{\rm tanh}\, p_2~+~b}\over{c~{\rm tanh}\, p_2~+~a}}\,,
\\\label{tanhs}
\\\nonumber-{\rm tanh}\, p_2= {{d~{\rm tanh}\, p_1~+~c}\over{b~{\rm tanh}\, p_1~+~a}}\,.
\end{eqnarray}
Plugging the latter equation into the former gives the quadratic
$$0=p_1^2(ab-cd)+p_1(a^2+b^2-c^2-d^2)+(ab-cd)$$
with discriminant
$$\begin{aligned}
\Delta&=(a^2+b^2-c^2-d^2)^2-4(ab-cd)^2\\
&=(a-b-c+d)(a-b+c-d)(a+b-c-d)(a+b+c+d)\\
&>0\\
\end{aligned}$$
since each factor is positive, as noted above since $L_1$ and $L_2$ are not ultraparallel.

\medskip

To see that this critical point is in fact a minimizer, we apply the second derivative
test for a function of two independent variables.  To this end, we have
$$\begin{aligned}
{{\partial^2{\rm cosh}\,d}\over{\partial p_1^2}}&={{\partial^2{\rm cosh}\,d}\over{\partial p_2^2}}={\rm cosh}\, d >0,\\\\
{{\partial^2{\rm cosh}\,d}\over{\partial p_1\partial p_2}}&={\rm sinh}\, p_2\{b~{\rm cosh}\,p_1~+~a~{\rm sinh}\,p_1\} \\
&+~{\rm cosh}\, p_2\{d~{\rm cosh}\,p_1~+~c~{\rm sinh}\,p_1\} ,\\
\end{aligned}$$
whence

$$\begin{aligned}
\biggl(&{{\partial^2{\rm cosh}\,d}\over{\partial p_1^2}}\biggr)\biggl({{\partial^2{\rm cosh}\,d}\over{\partial p_2^2}}\biggr)-\biggl({{\partial^2{\rm cosh}\,d}\over{\partial p_1\partial p_2}}\biggr)^2\\\\
&=\biggl[{\rm cosh}\,p_1(a~{\rm cosh}\,p_2~+~c~{\rm sinh}\, p_2)~+~{\rm sinh}\,p_1(b~{\rm cosh}\,p_2~+~d~{\rm sinh}\, p_2)\biggr]^2\\
&-\biggl[{\rm sinhh}\,p_2(b~{\rm cosh}\,p_1~+~a~{\rm sinh}\, p_1)~+~{\rm cosh}\,p_2(d~{\rm cosh}\,p_1~+~c~{\rm sinh}\, p_1)\biggr]^2\\
&= \biggl[{\rm cosh}\,p_1((a+d)~{\rm cosh}\,p_2~+~(b+c)~{\rm sinh}\, p_2)~\\
&\hskip 2cm+~{\rm sinh}\,p_1((b+c)~{\rm cosh}\,p_2~+~(a+d)~{\rm sinh}\, p_2)\biggr]\\
&\times\biggl[{\rm cosh}\,p_1((a-d)~{\rm cosh}\,p_2~-~(b-c)~{\rm sinh}\, p_2)~\\
&\hskip 2cm+~{\rm sinh}\,p_1((b-c)~{\rm cosh}\,p_2~-~(a-d)~{\rm sinh}\, p_2)\biggr]\\
&=\biggl[(a+d){\rm cosh}(p_1+p_2)+(b+c){\rm sinh}(p_1+p_2)\biggr]~\\
&\times\biggl[(a-d){\rm cosh}(p_1-p_2)+(b-c){\rm sinh}(p_1-p_2)\biggr]\\
&={\rm cosh}(p_1+p_2)~{\rm cosh}(p_1-p_2)\\
&\times\biggl [(a^2-d^2)+{{(b^2-c^2)^2}\over{a^2-d^2}}-(b+c)^2{{a-d}\over{a+d}}-(b-c)^2{{a+d}\over{a-d}}\biggr]\\
&={{{\rm cosh}(p_1+p_2)~{\rm cosh}(p_1-p_2)}\over{a^2-d^2}}~\Delta\\
&>0,
\end{aligned}$$
since $\Delta>0$ and $a^2>d^2$ by the remarks that began this proof, where we have used the
equalities (\ref{tanhs}) to write ${\rm sinh}(p_1\pm p_2)={{b\pm c}\over{a\pm d}}{\rm cosh}(p_1\pm p_2)$.

\bigskip

This completes the proof that the critical point is a minimizer, and since ${\rm cosh}\, d$ is convex for $d>0$, this guarantees a unique minimizer.  
It follows without difficulty from the Super Law of Cosines Corollary \ref{LOC} that this minimizer is indeed orthogonal to $L_1$ and $L_2$.
\end{proof}

The proof provides a formula for the distance between the two geodesics, which is omitted and  is far more complicated than the the classical case
of two parallel non-ultraparallel geodesics.  Just as for the Law of Cosines, the same result encodes substantial fermionic corrections.

\end{document}